\documentclass [12pt,a4paper]{article}
\usepackage[T1]{fontenc}
\usepackage[ansinew]{inputenc}
\usepackage{amssymb}
\usepackage{amsmath}
\usepackage{graphicx}
\usepackage{amsthm}
\usepackage{cite}

\newtheorem{theorem}{\textbf{Theorem}}[section]
\newtheorem{lemma}{\textbf{Lemma}}[section]
\newtheorem{proposition}{\textbf{Proposition}}[section]
\newtheorem{corollary}{\textbf{Corollary}}[section]
\newtheorem{remark}{\textbf{Remark}}[section]
\newtheorem{definition}{\textbf{Definition}}[section]

\allowdisplaybreaks[4]

\def\be{\begin{equation}}
\def\ee{\end{equation}}
\def\bea{\begin{eqnarray}}
\def\eea{\end{eqnarray}}
\def\bt{\begin{theorem}}
\def\et{\end{theorem}}
\def\bl{\begin{lemma}}
\def\el{\end{lemma}}
\def\br{\begin{remark}}
\def\er{\end{remark}}
\def\bp{\begin{proposition}}
\def\ep{\end{proposition}}
\def\bc{\begin{corollary}}
\def\ec{\end{corollary}}
\def\bd{\begin{definition}}
\def\ed{\end{definition}}

\def\s0t{\sup _{0\leq \tau\leq t}}
\def\C0T{C([0,T];\,}

\def\DAS{D( A^{\frac{S}{2}})}
\def\DAS1{D( A^{ \frac{S+1}{2}})}

 \def\non{\nonumber }

\def \no#1#2#3 {{\bf #1} (#3), #2.}
\def \eds#1#2#3 {#1, #2, #3.}

\begin{document}

\title{Long-time Behavior for Nonlinear Hydrodynamic System Modeling the Nematic Liquid Crystal Flows }

\author{{\sc Hao Wu}\\ School of
Mathematical
Sciences, Fudan University\\Shanghai 200433, China\\
haowufd@yahoo.com}

\date{\today}

\maketitle


\begin{abstract}
We study a simplified system of the original Ericksen-Leslie
equations for the flow of nematic liquid crystals. This is a coupled
non-parabolic dissipative dynamic system. We show the convergence of
global classical solutions to single steady states as time goes to
infinity (uniqueness of asymptotic limit) by using the \L
ojasiewicz--Simon approach. Moreover, we provide an estimate on the
convergence rate. Finally, we discuss some possible extensions of
the results to certain generalized problems with changing density or
free-slip boundary condition.

\end{abstract}

\section{Introduction}
We consider the following hydrodynamical model for the flow of
nematic liquid crystals (cf. \cite{lin1,LL95})
 \bea
 v_t+v\cdot\nabla v-\nu \Delta v+\nabla P&=&-\lambda
 \nabla\cdot(\nabla d\odot\nabla d),\label{1}\\
 \nabla \cdot v &=& 0,\label{2}\\
 d_t+v\cdot\nabla d&=&\gamma(\Delta d-f(d)),\label{3}
 \eea
in $\Omega \times\mathbb{R}^+$, where $\Omega \subset \mathbb{R}^n\
(n=2,3)$ is a bounded domain with smooth boundary $\Gamma$. Here,
$v$ is the velocity field of the flow and $d$ represents the
averaged macroscopic/continuum molecular orientations in
$\mathbb{R}^n \ (n=2,3)$. $P(x,t)$ is a scalar function representing
the pressure (including both the hydrostatic and the induced elastic
part from the orientation field). The positive constants $\nu,
\lambda$ and $\gamma$ stand for viscosity, the competition between
kinetic energy and potential energy, and macroscopic elastic
relaxation time (Debroah number) for the molecular orientation
field. We assume that $f(d)=\nabla F(d)$ for some smooth bounded
function $F:\mathbb{R}^n\rightarrow \mathbb{R} $. $\nabla d\odot
\nabla d$ denotes the $n\times n$ matrix whose $(i,j)$-th entry is
given by $\nabla_i d\cdot \nabla_j d$, for $1\leq i,j\leq n$.

In this paper we deal with the system \eqref{1}--\eqref{3} subject
to the initial conditions
 \be
 v|_{t=0}=v_0(x) \ \ \text{with}\ \nabla\cdot v_0=0,\quad
 d|_{t=0}=d_0(x),\qquad \text{for}\ x\in \Omega,\label{4}
 \ee
 and the Dirichlet boundary conditions:
 \be
 v(x,t)=0,\quad d(x,t)=d_0(x),\qquad \text{for}\ (x, t)\in \Gamma\times
 \mathbb{R}^+.
 \label{5}
 \ee

 In
\cite{lin1}, the author proposed equations \eqref{1}--\eqref{3} as a
simplified system of the original Ericksen--Leslie system (cf.
\cite{E1,Le}). By Ericksen--Leslie's hydrodynamical theory of the
liquid crystal, the (simplified) system describing the orientation
as well as the macroscopic motion reads as follows (here we assume
the density to be constant)
 \bea
 v_t+v\cdot\nabla v-\nu \Delta v+\nabla P&=&-\lambda
 \nabla\cdot(\nabla d\odot\nabla d),\label{a1}\\
 \nabla \cdot v &=& 0,\label{a2}\\
 d_t+v\cdot\nabla d&=&\gamma(\Delta d+|\nabla d|^2 d),\quad |d|=1.\label{a3}
 \eea
In order to avoid the gradient nonlinearly in \eqref{a3}, usually
one uses the Ginzburg--Landau approximation to relax the constraint
$|d|=1$. The corresponding approximate energy is
$$\int_\Omega \frac12|\nabla d|^2+\frac{1}{4\eta^2}(|d|^2-1)^2 dx,
$$
where $\eta$ is a positive constant. Then we arrive at the
approximation system \eqref{1}--\eqref{3}, where  \be
f(d)=\frac{1}{\eta^2}(|d|^2-1)d \label{FF}
 \ee with its
antiderivative  \be F(d)=\frac{1}{4\eta^2}(|d|^2-1)^2. \ee

The Ericksen--Leslie system is well suited for describing many
special flows for the materials, especially for those with small
molecules, and is wildly accepted in the engineering and
mathematical communities studying liquid crystals. System
\eqref{1}--\eqref{3} can be possibly viewed as the simplest
mathematical model, which keeps the most important mathematical
structure as well as most of the essential difficulties of the
original Ericksen--Leslie system (cf. \cite{LL95}). System
\eqref{1}--\eqref{3} with Dirichlet boundary conditions has been
studied in a series of work not only theoretically (cf.
\cite{LL95,LL96}) but also numerically (cf. \cite{Liu1,Liu2}). In
particular, in \cite{LL95}, the authors proved the existence theorem
for the weak solutions of system \eqref{1}--\eqref{5} by a modified
Galerkin scheme. After that they also obtained the global existence
and uniqueness of classical solutions to the same system for $ n =
2$ or $n = 3$ with large viscosity assumption. Moreover, a
preliminary analysis of the asymptotic behavior of global classical
solution was also given in \cite{LL95}. In the final remark of
\cite{LL95}, a natural question on the uniqueness of the asymptotic
limit was raised. This is just the main goal of the present paper.
For the sake of simplicity, in the following text, we always treat
the nonlinearity $f$ of form \eqref{FF}. However, it is not
difficult to verify that our results holds true for more general
nonlinearities which is analytic and with proper growth and
dissipation assumptions.

In this paper, we show the convergence to equilibrium of global
classical solutions to system \eqref{1}--\eqref{5}. Namely, we
obtain the following results:

 \bt
 \label{main2d}
 When $n=2$, for any $v_0\in H_0^1(\Omega)$ with $\nabla\cdot v_0=0$ and $d_0\in H^2(\Omega)$, the
 unique classical solution to problem \eqref{1}--\eqref{5} has the
following property
 \be \lim_{t\rightarrow +\infty}
 (\|v(t)\|_{H^1}+\|d(t)-d_\infty\|_{H^2})=0,\label{cgce}
 \ee
 where $d_\infty$ is a solution to the following nonlinear elliptic boundary value
  problem:
  \be \left\{\begin{array}{l}  - \Delta d_\infty + f(d_\infty)=0,\quad x\in \Omega, \\
   d_\infty=d_0(x), \;\;\; x\in \Gamma.\\
   \end{array}
   \label{staa}
  \right.
 \ee
 Moreover, there exists a positive constant $C$ depending on
 $v_0,d_0, \Omega, d_\infty$,
 such that
 \be
 \|v(t)\|_{H^1}+\|d(t)-d_\infty\|_{H^2}\leq C(1+t)^{-\frac{\theta}{(1-2\theta)}}, \quad \forall\ t \geq
 0,\label{rate}
 \ee
 with $\theta \in (0,1/2)$ being the same constant as in the \L
 ojasiewicz--Simon inequality (see Lemma \ref{ls} below).
 \et

When the spacial dimension is three, we deal with two cases. The
first result is concerning the large viscosity case, we have

\bt\label{main3da}
 When $n=3$, for any $v_0\in H_0^1(\Omega)$ with $\nabla\cdot v_0=0$,
 $ d_0\in H^2(\Omega)$ under the large viscosity assumption $\nu\geq
 \nu_0(\lambda, \gamma, v_0,d_0)$, the unique global classical solution of problem
 \eqref{1}--\eqref{5}
 enjoys the same properties as in Theorem \ref{main2d}.
 \et

 The second one is  a "stability" result for the near equilibrium
 initial data in the three dimensional case.

 \bt\label{main3d}
 When $n=3$, let $d^*\in H^2(\Omega)$ be an absolute minimizer of the functional
 $$ E(d)=\frac12\|\nabla d\|^2+ \int_\Omega F(d)dx$$
 in the sense that $E(d^*)\leq E(d)$ whenever $d=d^*=d_0(x)$ on $\Gamma$.
 There is a constant $\sigma$ which may depend on $\lambda, \gamma,
 \nu$ and $v_0,d_0$, such that if
 $\|v_0\|_{H^1}+\|d_0-d^*\|_{H^2}<\sigma$, then the problem
 \eqref{1}--\eqref{5} admits a unique global classical solution
 enjoying the same properties as in Theorem \ref{main2d}.
 \et
\br Theorem \ref{main3d} implies that if the initial data is
sufficiently close to an absolute minimizer of functional $E$, then
there exists a global solution and the solution will converge to an
equilibrium which may not necessarily be the original minimizer.
This is because the set of equilibria might be a continuum. Theorem
\ref{main3d} gives the "uniqueness" of asymptotic limit of  the
global solution to problem \eqref{1}--\eqref{5}. This improves the
result stated in \cite[Theorem C]{LL95}, in which only sequence
convergence for
 director field $d$ was obtained.
 \er

The problem about uniqueness of asymptotic limit for nonlinear
evolution equations, namely whether the global solution will
converge to an equilibrium as time tends to infinity, has attracted
a lot of interests of mathematicians. If the space dimension $n\geq
2$, it is known that the structure of the set of equilibria can be
nontrivial and may
   form a continuum for certain physically reasonable nonlinearities. The reader is
referred, for instance, to \cite[Rem. 2.3.13]{HA91}, where the
following two-dimensional equation $-\Delta u+u^3-\lambda u=0$,
$\lambda>0$, endowed with a standard Dirichlet homogeneous boundary
condition, is considered. And we note that, for the vector
functions, the situations may be even more complicated. If this is
the case, it is highly nontrivial to decide whether or not a given
bounded trajectory converges to a single steady state. In 1983, L.
Simon \cite{S83} made a breakthrough that for a semilinear parabolic
equation with
  a nonlinearity $f(x,u)$ being analytic in the unknown function $u$, its bounded global solution would
  converge to an equilibrium  as $t\rightarrow \infty$.  Simon's idea relies on a generalization of the
  \L ojasiewicz inequality (see \cite{L1,L2}) for analytic functions defined in finite
  dimensional space $\mathbb{R}^m$. Since then, his original approach has been simplified and applied
  to prove convergence results for many evolution equations (see e.g.,
\cite{LD,J981,HJ99,HJ01,ET01,WGZ1,WGZ3,WGZ4,Hu,W07,WZ1} and the
references cited therein). For our problem \eqref{1}--\eqref{5}, in
order to apply the \L ojasiewicz--Simon approach to prove the
convergence result, we need to introduce a suitable
{\L}ojasiewicz--Simon type inequality for vector functions with
nonhomogeneous Dirichlet boundary condition (cf. Lemma \ref{ls}).

As far as the convergence rate is concerned, it is known that an
estimate in certain (lower order) norm can usually be obtained
directly from the \L ojasiewicz--Simon approach (see, e.g.,
\cite{Z04,HJ01}). Then, one straightforward way to get estimates in
higher order norms is using interpolation inequalities (cf.
\cite{HJ01}) and, consequently, the decay exponent deteriorates. We
shall show that by using suitable energy estimates and constructing
proper differential inequalities, it is possible to obtain the same
estimates on convergence rate in both higher and lower order norms.
 Our approach
in some sense improves the previous results in the literature (see,
for instance, \cite{HJ01,Z04}) and it can
 apply to many other problems (cf. \cite{WGZ1,WGZ3,WGZ4,W07}).

 The remaining part of this paper is organized as follows. In Section 2,
 we introduce the functional setting, some preliminary results
  as well as some technical lemmas.  Section 3 is devoted to the two dimensional case. We prove
  the convergence of global solutions to
single steady states as time goes to infinity and obtain an estimate
on convergence rate. In Section 4, we consider the three dimensional
case. The same convergence result was proved for two subcases, in
which the global existence of classical solutions can be obtained.
In the final Section 5, we discuss some possible extensions of our
results to certain generalized problems with changing density or
free-slip boundary conditions.

\section{Preliminaries}
\setcounter{equation}{0} First, we introduce the function spaces we
shall work on (cf. \cite{LL95,Te}):
 \bea && H^1_0(\Omega)=\text{the\ closure\ of}\
C_0^\infty(\Omega,\mathbb{R}^n)\ \text{in\ the\ norm}\
\left(\int_\Omega|\nabla v|^2dx\right)^\frac12,\non\\
 && H^{-1}(\Omega)
= \text{the\ dual\ of} \ H^1_0(\Omega),\non\\
 &&  H^2(\Omega)= \{v\in L^2(\Omega, \mathbb{R}^n)\ |\
v_{x_i},v_{x_ix_j}\in L^2(\Omega, \mathbb{R}^n),\ 1\leq i,j\leq
n\},\non\\
 && \mathcal{V}=C_0^\infty(\Omega,\mathbb{R}^n)\cap \{v: \nabla\cdot v=0\},\non\\
 && H=\ \text{the closure of}\ \mathcal{V} \text{\ in\ } L^2(\Omega,
\mathbb{R}^n),\non\\
 && V=\ \text{the closure of\ } \mathcal{V} \text{\ in\ } \ H^1_0(\Omega),\non\\
 && V'=\text{the\ dual of\ } V.\non
 \eea

Global existence and uniqueness of classical solution to system
\eqref{1}--\eqref{5} has been proven in \cite[Theorem B]{LL95}. More
precisely, we have
 \bp\label{glo1}
 Problem \eqref{1}--\eqref{5} admits a unique global classical
 solution $(v,d)$ provided that $v_0\in H_0^1(\Omega), d_0\in H^2(\Omega)$ either
 $n=2$ or $n=3$ with the large viscosity assumption $\nu\geq
 \nu_0(\lambda, \gamma, v_0,d_0)$.
 \ep

For any classical solution $(v,d)\in \Omega\times [0,T]=Q_T$ $
(0\leq T\leq +\infty)$ of problem \eqref{1}--\eqref{5}, we consider
the functional
 \be \mathcal{E}(t)=\frac{1}{2}\|v(t)\|^2+\frac{\lambda}{2}\|\nabla
 d(t)\|^2+\frac{\lambda}{2}\int_\Omega F(d(t))dx.\label{Ly}
 \ee
 It has been shown in \cite{LL95} that our system \eqref{1}--\eqref{5} has the following \textit{basic
energy law}, which can be viewed as a direct consequence of the
balance laws of the linear momentum \eqref{1} and angular momentum
\eqref{3}:
 \be
 \frac{d}{dt}\mathcal{E}(t)+\nu\|\nabla v(t)\|^2+\lambda\gamma\|\Delta
 d(t)-f(d(t))\|^2=0,
 \quad 0\leq t\leq T.\label{lya}
 \ee
\eqref{lya} reflects the energy dissipation property of the flow of
liquid crystals. Moreover, one can verify that $\mathcal{E}(t)$
serves as a Lyapunov functional for problem \eqref{1}--\eqref{5}.

Next, we look at the following elliptic boundary value problem
\be \left\{\begin{array}{l}  - \Delta d + f(d)=0,\quad x\in \Omega, \\
   d=d_0(x), \;\;\; x\in \Gamma.\\
   \end{array}
   \label{staaq}
  \right.
 \ee
Denote
 \be
 E(d)=\frac12\|\nabla d\|^2 + \int_\Omega F(d)dx.\label{EDD}
 \ee
 It is not difficult to see that the solution to \eqref{staaq} is a
 critical point of $E(d)$, and conversely, the
 critical point of $E(d)$ is a solution to \eqref{staaq} (cf. \cite{WZ1,WGZ1} and references cited
 therein). Besides, regularity of the solution to \eqref{staaq} has been shown in \cite{LL95}
 such that $d$ is smooth on $\Omega$ provided $d_0$ is smooth on $\Gamma$.

  As mentioned in Introduction, in order to apply the \L ojasiewicz--Simon approach to prove the
 convergence to equilibrium, we have to introduce a suitable \L
 ojasiewicz--Simon type inequality related to our present problem. In
 particular, we have
 \bl\label{ls}{\rm [\L ojasiewicz--Simon Type Inequality]}
 Let $\psi$ be a critical point of $E(d)$. There exist constants
 $\theta\in(0,\frac12)$ and $\beta>0$ depending on $\psi$ such that
 for any $d\in H^1(\Omega)$ satisfying $d|_\Gamma=d_0(x)$ and $\|d-\psi\|_{H^1}<\beta$,
 there holds
 \be
 \|-\Delta d+f(d)\|_{H^{-1}}\geq
 |E(d)-E(\psi)|^{1-\theta}.
 \ee
 \el

 \br
 The above lemma can be viewed as an extended version of Simon's result \cite{S83} for scalar function under the use of $L^2$-norm.
 We can refer to \cite[Chapter 2, Theorem 5.2]{Hu}, in which the case
 for vectors subject to homogeneous Dirichlet boundary condition was considered. Here
 we observe that, our present (nontrivial) boundary data for  director field $d$ does not depend
 on time. As a result, every solution to the corresponding stationary problem \eqref{staaq}, which
 is a critical point of $E(d)$ satisfies the same boundary condition as the solution to the evolution problem. Therefore, we
 only have to derive a \L ojasiewicz--Simon type inequality for functions $d$, which satisfy $d|_\Gamma=d_0(x)$ and
 fall into a properly small neighborhood of certain but arbitrary critical point of
 $E(d)$. In this case, it is always true that the difference $\tilde{d}=d-\psi\in H_0^1(\Omega)$.
 Keeping this fact in mind, we are able to prove the present lemma following the steps in \cite[Chapter 2, Theorem
 5.2]{Hu} or \cite{HJ99}. Hence,  the details are omitted here. We could also refer to a related case treated in
 \cite{WGZ3}, that a  nonhomogeneous (time-dependent Dirichlet) boundary condition was removed by
 a proper variable transformation (cf. also \cite{WZ1} where the boundary condition contains a nonzero constant).
  \er

In the following text, we will use the regularity result for Stokes
problem (cf. \cite{Te1})
  \bl \label{S}
  Denote the Stokes operator by $S$, which is a unbounded operator
  in $H$ of domain $H^2(\Omega)\cap V$:
  $$Su=-\Delta u +\nabla \pi \in \ H, \quad \forall u \in H^2(\Omega)\cap
  V.$$
  Then there exists a constant $C$ such that for any $u\in
  H^2(\Omega)\cap V$,
  $$\|u\|_{H^2}+\|\pi\|_{H^1\setminus\mathbb{R} }\leq C\|S u\|.$$
  \el

 Before ending this section, we introduce the following lemma which is useful in the study of large time behavior
 of solutions to evolution problems. We will apply it to obtain
  uniform (higher order) estimates of the solution and decay of the energy dissipations
  of system \eqref{1}--\eqref{5}.

 \bl \label{SZ} {\rm \cite[Lemma 6.2.1]{Z04}}
 Let $T$ be given with $0<T\leq +\infty$. Suppose that $y(t)$ and $h(t)$ are nonnegative continuous functions defined on
 $[0,T]$ and satisfy the following conditions:
 $$
 \frac{dy}{dt}\leq c_1 y^2+ c_2 +h(t),\quad \text{with}\ \
 \int_0^T y(t) dt\leq c_3,\quad \int_0^T h(t)dt\leq c_4,
 $$
 where $c_i(i=1,2,3,4)$ are given nonnegative constants. Then for
 any $r\in (0,T)$, the following estimates holds:
 $$
 y(t+r)\leq \left(\frac{c_3}{r}+c_2r+c_4\right)e^{c_1c_3},\quad
 \forall\ t\in[0,T-r].
 $$
 Furthermore, if $T=+\infty$, then
 $$\lim_{t\rightarrow +\infty} y(t)=0.$$
 \el

\section{Convergence to Equilibrium for Two Dimensional Case }
\setcounter{equation}{0}

In this section, we prove the convergence of global solutions to
single steady states as time tends to infinity for 2-D case. Since
parameters $\lambda, \gamma,\nu$ do not play crucial role in the 2-D
case, we set $\lambda=\gamma=\nu=1$ in this section for the sake of
simplicity.

When the space dimension equals to two, an important property for
the global solution to problem \eqref{1}--\eqref{5} is the following
high order energy law, which played a crucial role in the proof of
global existence result in \cite{LL95}. Denote
 \be
 A(t)=\|\nabla v(t)\|^2+ \|\Delta d(t)-f(d(t))\|^2.\label{A}
 \ee

 Then we have
 \bl \label{he2d} (cf. \cite[(4.9)]{LL95}) In 2-D case, the following inequality holds for the
 classical solution $(v,d)$ to problem \eqref{1}--\eqref{5}
 \be \frac{d}{dt}A(t)+(\|\Delta v\|^2+\|\nabla(\Delta d-f(d))\|^2)\leq C(A^2(t)+A(t)), \quad \forall
 \
 t\geq 0,\label{he}
 \ee
 where $C$ is a constant depending on $f, \Omega, \|v_0\|, \|d_0 \|_{H^1(\Omega)}$.
 \el

\subsection{Convergence to Equilibrium}

Based on the high order energy law \eqref{he}, we are able to show
the convergence of the velocity field $v$ first.
 \bl \label{vcon}
 For any $t\geq 0$, the following uniform estimate holds
 \be \|v(t)\|_{H^1}+ \|d(t)\|_{H^2}\leq C, \label{ubdd}
 \ee
 where $C$ is a constant depending on $f, \Omega, \|v_0\|_{H^1}, \|d_0
 \|_{H^2(\Omega)}$.  Furthermore,
 \be \lim_{t\rightarrow +\infty} (\|v(t)\|_{H^1}+ \|-\Delta d(t)+f(d(t))\|)=0. \label{vcon1}\ee
 \el
\begin{proof}
It follows from the basic energy law \eqref{Ly} that
 \be \mathcal{E}(t)+ \int_0^t A(\tau) d\tau =
 \mathcal{E}(0)<\infty,\quad \forall\ t\geq 0.\ee
  By the Young inequality $a^2\leq \frac12 a^4+ \frac12$, we can see
  that $\mathcal{E}(t)$ is bounded from below by a constant which is only dependent of $|\Omega|$. As a result,
  \be \int_0^\infty A(t)dt\leq \mathcal{E}(0)< +\infty,\label{ae}
  \ee
  and
  \be \mathcal{E}(t)\leq \mathcal{E}(0),\quad  \forall\ t\geq 0.\label{eebd}\ee
  \eqref{eebd} implies the uniform estimate
  \be \|v(t)\|+\|d(t)\|_{H^1}\leq C,\quad  \forall\ t\geq 0.\label{bd1}\ee
 Furthermore, \eqref{ae} together with Lemma \ref{he2d}  and Lemma \ref{SZ} yields
  that
  \be \lim_{t\rightarrow +\infty} \left(\|\nabla v(t)\|+ \|\Delta
  d(t)-f(d(t))\|\right)=0. \label{con1}
  \ee
  By the Poincar\'e inequality, we prove the conclusion
  \eqref{vcon1}. Concerning the uniform bound \eqref{ubdd}, we take $r=1$ in Lemma \ref{SZ}
  to get
  \be
  \|\nabla v(t)\|+\|-\Delta d(t)+f(d(t))\|\leq C,\quad \forall \ t\geq 1,
  \label{bd2}
  \ee
 where $C$ does not depend on $t$. On the other hand, for any $t\in
 [0,1]$, it follows from \eqref{he} and the fact $\int_0^1
 A(t)dt\leq C$ that
 \be \sup_{0\leq t\leq 1} A(t) \leq e^{\int_0^1 A(t)dt}A(0)+C\leq C. \label{bd4}\ee
 Besides, from the continuous embedding $H^1\hookrightarrow L^p (1\leq p<\infty)$ and \eqref{bd1} we have
 \be \|\Delta d\|\leq \|-\Delta d+f(d)\|+\|f(d)\|\leq \|-\Delta d+f(d)\|+
 C\left(1+\|d\|_{L^6}^3\right)\leq C. \label{bd3}
 \ee
 Now we can conclude \eqref{ubdd} from
 \eqref{bd2}--\eqref{bd3}. The proof is complete.
\end{proof}

 Let $\mathcal{S}$ be the set
 $$ \mathcal{S}=\{(0,u)\ | \ -\Delta u+ f(u)=0,\
 \text{in} \ \Omega, \ u|_\Gamma=d_0(x) \}.$$
The $\omega$-limit set of $(v_0,d_0)\in V\times H^2(\Omega)\subset
L^2(\Omega) \times H^1(\Omega)$ is defined as follows:
 \bea
 \omega((v_0,d_0)) &= &\{ (v_\infty(x),d_\infty(x)) \mid\ \text{there
 \ exists\ } \{t_n\}\nearrow \infty  \text{\ such\ that\ } \non\\&&
 (v(x,t_n),d(x,t_n)) \rightarrow (v_\infty(x),d_\infty(x))\
 \text{in}\ L^2 \times H^1,\ \text{as}\ t_n\rightarrow +\infty
 \}.\non
 \eea

 We infer from Lemma  \ref{vcon} that
 \begin{proposition} \label{lim}  $\omega((v_0,d_0))$ is a nonempty bounded subset in $H^1(\Omega)\times
 H^2(\Omega)$. Besides, all asymptotic
 limiting  points $(v_\infty, d_\infty)$ of problem \eqref{1}--\eqref{5} belong to
 $\mathcal{S}$. In other words, $\omega((v_0,d_0))\subset
 \mathcal{S}$.
 \end{proposition}

In what follows, we prove the convergence for director field $d$.
For any initial datum $(v_0,d_0)\in V\times H^2(\Omega)$, it follows
from Lemma \ref{vcon} that $ \|d\|_{H^2}$ is uniformly bounded.
 Since the embedding $H^2\hookrightarrow H^1$ is compact,
 there is an increasing unbounded sequence
$\{t_n\}_{n\in\mathbb{N}}$ and a function $d_\infty$ such that
   \be \lim_{t_n\rightarrow +\infty} \|d(t_n)-d_\infty\|_{H^1}
   =0. \label{secon}
   \ee
In particular, Proposition \ref{lim} implies that $d_\infty$
satisfies the equation
  \be -\Delta d_\infty+f(d_\infty)=0,\quad x\in \Omega,\quad
  d_\infty|_\Gamma=d_0.\label{sta}
  \ee

We prove the convergence result following a simple argument
introduced in \cite{J981}, in which the key observation is that
after a certain time $t_0$,  $d(t)$ will fall into a certain small
neighborhood of $d_\infty$ and stay there forever.

From the basic energy law \eqref{Ly}, we can see that
$\mathcal{E}(t)$ is decreasing on $[0,\infty)$, and it has a finite
limit as time goes to infinity because it is bounded from below.
Therefore, it follows from \eqref{secon} that
 \be \lim_{t_n\rightarrow +\infty} \mathcal{E}(t_n)=E(d_\infty).
 \ee On the other hand,  we can infer from \eqref{Ly} that
$\mathcal{E}(t) \geq E(d_\infty)$, for all $t>0$, and the equal sign
holds if and only if, for all $t>0$, $v=0$ and $d$ solves problem
(\ref{sta}).

We now consider all possibilities.

\noindent \textbf{ Case 1}. If there is a $t_0>0$ such that at this
time
   $\mathcal{E}(t_0)=E(d_\infty)$, then for all $t>t_0$, we deduce from
   (\ref{Ly}) that
   \be \|\nabla v\|\equiv 0,\quad  \|-\Delta d+f(d)\|\equiv 0. \label{muee}\ee
It follows from \eqref{3}, \eqref{muee} and the Sobolev embedding
Theorem that for $t>t_0$
 \be
 0\leq \|d_t\|\leq \|v\cdot\nabla d\|+ \|-\Delta d+f(d)\|\leq
 \|v\|_{L^4}\|\nabla d\|_{L^4}\leq C\|\nabla v\|=0.
 \ee
Namely, $d$ is independent of time for all $t>t_0$.  Due to
\eqref{secon}, we have $d(t)\equiv d_\infty$ for $t>t_0$.

\noindent \textbf{Case 2}. For all $t>0$,
$\mathcal{E}(t)>E(d_\infty)$. First we assume that the following
claim holds true.
 \begin{proposition} \label{neighb} There is a $t_0>0$ that for all $t\geq t_0$,
$\|d(t)-d_\infty\|_{H^1}<\beta$. Namely, for all $t\geq t_0$, $d(t)$
satisfies the condition in Lemma \ref{ls}.
 \end{proposition}

  In this
case, it follows from Lemma \ref{ls} that
 \be
 |E(d)-E(d_\infty)|^{1-\theta}\leq \|-\Delta d+f(d)\|_{H^{-1}}\leq \|-\Delta
 d+f(d)\|,\quad \forall\ t\geq t_0.
 \ee
The fact $\theta\in(0,\frac12)$ implies $0<1-\theta<1$,
$2(1-\theta)>1$. As a consequence,
  $$ \|v\|^{2(1-\theta)}=\|v\|^{2(1-\theta)-1}\|v\|\leq C\|v\|.$$
Then we infer from the basic inequality
$$ (a+b)^{1-\theta}\leq a^{1-\theta}+b^{1-\theta}, \qquad \forall \ a,b\geq 0$$
that
 \bea (\mathcal{E}(t)-E(d_\infty))^{1-\theta}&\leq&
 \left(\frac12\|v\|^2+|E(d)-E(d_\infty)|\right)^{1-\theta}\non\\
 &\leq& \left(\frac12\|v\|^2+\|-\Delta
 d+f(d)\|^{\frac{1}{1-\theta}}\right)^{1-\theta}\non\\
 &\leq&\left(\frac12\right)^{1-\theta} \|v\|^{2(1-\theta)}+\|-\Delta
 d+f(d)\|\non\\
 &\leq& C \|v\|+\|-\Delta
 d+f(d)\|.
 \eea
Therefore, a direct calculation yields
 \bea
 -\frac{d}{dt}(\mathcal{E}(t)-E(d_\infty))^\theta&=&
 -\theta(\mathcal{E}(t)-E(d_\infty))^{\theta-1}\frac{d}{dt}\mathcal{E}(t)\non\\
   &\geq &\frac{C\theta(\|\nabla v\|+ \|-\Delta d+f(d)\|)^2}{C \|v\|+\|-\Delta
 d+f(d)\|}\non\\
 &\geq& C_1(\|\nabla v\|+ \|-\Delta d+f(d)\|), \quad \forall\ t\geq
 t_0,\label{ly2}
 \eea
where $C_1$ is a constant depending on $v_0,d_0,\Omega$.\\
Integrating from $t_0$ to $t$,   we get
 \bea && (\mathcal{E}(t)-E(d_\infty))^\theta+ C_1\int_{t_0}^t (\|\nabla v(\tau)\|+ \|-\Delta
 d(\tau)+f(d(\tau))\|)d\tau\non\\
 & \leq&  (\mathcal{E}(t_0)-E(d_\infty))^\theta<\infty,
 \quad \forall\  t\geq t_0.\label{int}
 \eea
Since $\mathcal{E}(t) - E(d_\infty) \geq 0$, we conclude that
 \be  \int_{t_0}^\infty (\|\nabla v(\tau)\|+ \|-\Delta
 d(\tau)+f(d(\tau))\|)d\tau<\infty.
 \ee
On the other hand, it follows from equation \eqref{3} that
 \bea \|d_t\|&\leq& \|v\cdot \nabla d\| + \|-\Delta d+f(d)\leq \|v\|_{L^4}\|\nabla d\|_{L^4}+\|-\Delta d+f(d)\|\non\\
 &\leq & C\|\nabla v\|+ \|-\Delta d+f(d)\|.\label{dt}
 \eea
Hence,
 \be \int_{t_0}^\infty \|d_t(\tau)\| d\tau<+\infty,
 \ee
which easily implies that as $t\rightarrow +\infty$, $d(x, t)$
converges in $L^2(\Omega)$. This and \eqref{secon} indicate that
 \be \lim_{t\rightarrow +\infty}
 \|d(t)-d_\infty\|=0.
 \ee
 Since $d(t)$ is uniformly bounded in $H^2(\Omega)$ (cf. \eqref{ubdd}), by interpolation we have
 \be \lim_{t\rightarrow +\infty}
 \|d(t)-d_\infty\|_{H^1}=0.\label{conh1}
 \ee On the other hand, uniform bound of $d$ in $H^2(\Omega)$ implies the weak convergence
  $$ d(t)\rightharpoonup d_\infty,\quad \text{in}\ H^2(\Omega).$$
 However, the decay property of the quantity $A(t)$ (cf. Lemma \ref{vcon}) could tell us
 more. Namely, we could get strong convergence of $d$ in $H^2$
 without using uniform estimates in higher order norm. To see this, we
 keep in mind that
 that
 \bea \|\Delta d- \Delta d_\infty\|&\leq& \| \Delta d- \Delta d_\infty
 -f(d)+f(d_\infty)\|+ \|f(d)-f(d_\infty)\|\non
 \\
 &\leq& \|\Delta
 d-f(d)\|+\|f'(\xi)\|_{L^4}\|d-d_\infty\|_{L^4}\non\\
 &\leq& \|\Delta
 d-f(d)\|+C\|d-d_\infty\|_{H^1}.\label{kkk}
 \eea
 The above estimate together with \eqref{vcon1} and \eqref{conh1} yields
 \be \lim_{t\rightarrow +\infty}
 \|d(t)-d_\infty\|_{H^2}=0.\label{conh2}
 \ee
To finish the proof, we will show that Proposition \ref{neighb}
always holds true for the global solution $d(t)$ to system
\eqref{1}--\eqref{5}. Define
 \be \bar{t}_n = \sup \{\ t > t_n | \ \|d(\cdot,s) - d_\infty\|_{H^1} < \beta ,\ \forall \ s\in [t_n ,
 t] \}. \label{tn}
 \ee
 It follows from (\ref{secon}) that for  any $ \varepsilon\in (0, \beta)$,  there exists  an integer $N$ such that when
$n \geq N$,
 \bea && \| d(\cdot,t_n) - d_\infty \|_{H^1} < \varepsilon, \label{sea}\\
 &&\frac{1}{C_1} (\mathcal{E}(t_n) - E(d_\infty))^\theta <
 \varepsilon. \label{seb}
 \eea
On the other hand, we can easily see that the orbit of $d$ is
continuous in $H^1$. This is because we already know from
\eqref{ubdd} that $d\in L^\infty(0,+\infty; H^2(\Omega))$. As a
consequence, $d\in L^2(t,t+1; H^2(\Omega))$ for any $t\geq 0$. The
basic energy law and \eqref{dt} imply $d_t\in L^2(t,t+1;
L^2(\Omega))$. Thus, $ d\in C([t,t+1]; H^1(\Omega))$, for any $t\geq
0$ (cf. \cite{Evans}). The continuity of the orbit of $d$ in $H^1$
and (\ref{sea}) yield that
$$\bar{t}_n >t_n,\quad \text{for\  all\ } \ n \geq N.$$

 Then there are two possibilities:\\
(i). If there exists $n_0\geq N$ such that $\bar{t}_{n_0}=+\infty$,
then from the previous discussions in Case 1 and Case 2, the
theorem is proved.\\
(ii). Otherwise, for all $n\geq N$, we have $t_n < \bar{t}_n <
+\infty$, and for all $t\in [t_n, \bar{t}_n]$,
$E(d_\infty)<\mathcal{E}(t)$. Then from (\ref{int}) with $t_0$ being
replaced by $t_n$, and $t$ being replaced by $\bar{t}_n$, we get
from \eqref{seb} that
 \be \int_{t_n}^{\bar{t}_n} (\|\nabla v(\tau)\|+ \|-\Delta
 d(\tau)+f(d(\tau))\|)d\tau< \varepsilon.\label{int1}
 \ee
Thus, it follows that (cf. \eqref{dt})
 \bea \| d(\bar{t}_n) - d_\infty\| &\leq&  \|d(t_n)
-d_\infty \| + \int_{t_n}^{\bar{t}_n} \| d_t(\tau) \| d\tau\non\\
  &\leq & \|d(t_n)-d_\infty \|+C\int_{t_n}^{\bar{t}_n} (\|\nabla v(\tau)\|+ \|-\Delta
 d(\tau)+f(d(\tau))\|)d\tau\non\\
 &<& C\varepsilon,
 \eea
which implies that $\lim_{n\rightarrow +\infty} \|d(\bar{t}_n)
-d_\infty \|=0.$ Since $d(t)$ is relatively compact in
$H^1(\Omega)$, there exists a subsequence of $\{d(\bar{t}_n)\}$,
still denoted  by $\{d(\bar{t}_n)\}$ converging to $d_\infty$ in
$H^1(\Omega)$, i.e., when $n$ is sufficiently large,
$$\| d(\bar{t}_n) - d_\infty\|_{H^1}\ < \beta$$
which contradicts the definition of $\bar{t}_n$ that $\|d(\cdot,
\bar{t}_n)-d_\infty\|_{H^1}=\beta$.

Summing up, we have considered all the possible cases and the conclusion \eqref{cgce} is proved.\\

\subsection{Convergence Rate}
In this part, we shall show the estimate on convergence rate
\eqref{rate}. This can be achieved in several steps.

\noindent \textbf{Step 1.} As has been shown in the literature (cf.
for instance, \cite{Z04,HJ01}), an estimate on the convergence rate
in certain lower order norm could be obtained directly from the \L
ojasiewicz--Simon approach. From Lemma \ref{ls} and \eqref{ly2}, we
have
 \bea
 \frac{d}{dt}(\mathcal{E}(t)-E(d_\infty))+ C_1(\mathcal{E}(t)-E(d_\infty))^{2(1-\theta)}\leq 0, \quad \forall\ t\geq
 t_0,\label{ly3}
 \eea
which implies
 \be \mathcal{E}(t)-E(d_\infty)\leq
 C(1+t)^{-\frac{1}{1-2\theta}}\quad \forall\ t\geq
 t_0. \ee
  Integrating \eqref{ly2} on
$(t,\infty)$, where $t\geq t_0$,  it follows from \eqref{dt} that
 \be \int_t^\infty \|d_t\| d\tau \leq\int_t^\infty(C\|\nabla v\|+
 \|-\Delta f+f(d)\|)d\tau\leq
 C(1+t)^{-\frac{\theta}{1-2\theta}}.\ee
By adjusting the constant $C$ properly, we obtain
 \be
    \|d(t)-d_\infty\|\leq C(1+t)^{-\frac{\theta}{1-2\theta}}, \quad t\geq 0.\label{rate1}
 \ee

\noindent \textbf{Step 2.} In Step 1, we only obtain the convergence
rate of $d$ (in $L^2$). Unlike for the temperature variable in some
phase-field systems (cf. \cite{WGZ1,WGZ3} and references cited
therein), although we have got some decay information for the
velocity field $v$ such that
 \be
\int_t^\infty\|\nabla v\|d\tau\leq
 C(1+t)^{-\frac{\theta}{1-2\theta}},\ee
 it is not easy to prove convergence rate of $v$
 directly. This is because now $v$ satisfies a Navier--Stokes type
 equation, which is much more complicated than the heat equation
 for the temperature variable in phase-field systems. As a
 result, one cannot easily obtain relation between $\|\nabla v\|$ and
 $v_t$ (in certain possible norm) from the equation itself. However, it is possible to
 achieve our goal by using the idea in \cite{WGZ1}, where we use
 higher order energy estimates and construct proper differential
 inequalities (cf. also \cite{W07,WGZ3,WGZ4}). Besides, in this way
 the convergence rate of $d$ in higher order norm can be proved
 simultaneously.

The steady state solution corresponding to problem
\eqref{1}--\eqref{5} satisfies the following system (cf.
\cite{LL95})
 \bea
v_\infty\cdot\nabla v_\infty-\nu \Delta v_\infty+\nabla P_\infty&=&-
 \nabla\cdot(\nabla d_\infty\odot\nabla d_\infty),\label{s1}\\
 \nabla \cdot v_\infty &=& 0,\label{s2}\\
 v_\infty\cdot\nabla d_\infty&=&\Delta
 d_\infty-f(d_\infty),\label{s3}\\
 v_\infty|_\Gamma=0, && d_\infty|_\Gamma=d_0(x).\label{s4}
 \eea

Lemma \ref{vcon} implies that the limiting point of system
\eqref{1}--\eqref{5} has the form $(0,d_\infty)\in \mathcal{S}$. As
a result, system \eqref{s1}--\eqref{s4} can be reduced to
 \bea
 \nabla P_\infty&=&-\nabla d_\infty\cdot \Delta d_\infty-\nabla\left(\frac{|\nabla d_\infty|^2}{2}\right),\label{s1a}\\
 -\Delta d_\infty+f(d_\infty)&=&0,\label{s2a}\\
 d_\infty|_\Gamma&=&d_0(x),\label{ss3}
 \eea
where in \eqref{s1a} we have used the fact that $$\nabla\cdot(\nabla
d_\infty\odot\nabla d_\infty)=\nabla\left(\frac{|\nabla
d_\infty|^2}{2}\right) + \nabla d_\infty\cdot\Delta d_\infty.$$
Subtracting the stationary problem \eqref{s1a}--\eqref{ss3} from the
evolution problem \eqref{1}--\eqref{5}, we get
 \bea
 && v_t+v\cdot\nabla v-\nu \Delta v+\nabla (P-P_\infty)
 + \nabla \left(\left(\frac{|\nabla d|^2}{2}\right)-\left(\frac{|\nabla
 d_\infty|^2}{2}\right)\right)\non\\
 &=&-\nabla d\cdot \Delta d +\nabla d_\infty\cdot \Delta d_\infty ,\label{11}\\
 && \qquad \ \ \nabla \cdot v = 0,\label{22}\\
 && d_t+v\cdot\nabla d=\Delta (d-d_\infty)
 -f(d)+f(d_\infty),\label{33}\\
 && (d-d_\infty)|_\Gamma=0.
 \eea
Multiplying \eqref{11} by $v$ and \eqref{33} by $-\Delta
d+f(d)=-\Delta (d-d_\infty)+f(d)-f(d_\infty)$ respectively,
integrating on $\Omega$, and adding the results together, we obtain
 \bea
  && \frac{d}{dt} \left(\frac12\|v\|^2+\frac12\|\nabla d-\nabla
  d_\infty\|^2+\int_\Omega F(d)-F(d_\infty) -
  f(d_\infty)(d-d_\infty)dx\right)\non\\
  && \ +\nu\|\nabla v\|^2+ \|\Delta d-f(d)\|^2\non\\
  &  =&  (v, \nabla d_\infty\cdot \Delta d_\infty) \non\\
  &  =&  (v, \nabla d_\infty\cdot (\Delta d_\infty-f(d_\infty)))+ (v\cdot\nabla d_\infty,
  -f(d_\infty))\non\\
  &=&0.
  \label{ra1}
  \eea
Multiplying \eqref{33} by $d-d_\infty$ and integrating in $\Omega$,
we have
 \be \frac12\frac{d}{dt}\|d-d_\infty\|^2+\|\nabla (d-d_\infty)\|^2=
 -(v\cdot\nabla d,d-d_\infty) -(f(d)-f(d_\infty), d-d_\infty):=I_1.\label{ra3}\ee
 The right hand side can be estimated as follows
 \bea
 |I_1|&\leq& \|v\|_{L^4}\|\nabla
 d\|_{L^4}\|d-d_\infty\|+\|f'(\xi)\|_{L^3}\|d-d_\infty\|^2_{L^3}\non\\
 &\leq& C\|\nabla v\|\|d-d_\infty\|+
 C(\|\nabla(d-d_\infty)\|^\frac13\|d-d_\infty\|^\frac23+\|d-d_\infty\|)^2\non\\
 &\leq& \varepsilon_1\|\nabla v\|^2+
 \frac12\|\nabla(d-d_\infty)\|^2+C\|d-d_\infty\|^2.\label{ra4}
 \eea
Multiplying \eqref{ra3} by $\alpha>0$ and adding the resultant to
\eqref{ra1}, using \eqref{ra4} we get
 \bea
  && \frac{d}{dt} \left(\frac12\|v\|^2+\frac12\|\nabla d-\nabla
  d_\infty\|^2+ \frac{\alpha}{2} \|d-d_\infty\|^2 +\int_\Omega F(d)dx-\int_\Omega F(d_\infty) dx\right.\non\\
  &&\ \ \left.-
  \int_\Omega f(d_\infty)(d-d_\infty)dx\right)+\left(\nu-\alpha\varepsilon_1\right)\|\nabla v\|^2
  + \|\Delta d-f(d)\|^2  +\frac{\alpha}{2}\|\nabla(d-d_\infty)\|^2  \non\\
  &\leq& C\alpha\|d-d_\infty\|^2.\label{ra5}
  \eea
On the other hand, by the Taylor's expansion, we have
 \be
F(d)=F(d_\infty) + f(d_\infty)(d-d_\infty) + f'(\xi)(d-d_\infty)^2,
\ee
where $\xi=a d+ (1-a)d_\infty$ with $a\in [0,1]$.\\
 Then we deduce that
 \bea & & \left\vert\int_\Omega F(d)dx-\int_\Omega
        F(d_\infty)dx+\int_\Omega f(d_\infty)d_\infty
         dx-\int_\Omega  f(d_\infty)d\  dx\right\vert\non\\
      &  =  & \left\vert \int_\Omega f'(\xi)(d-d_\infty)^2 dx\right\vert\non\\
      & \leq &
          \|f'(\xi)\|_{L^\infty}\|d-d_\infty\|^2\leq C_2\|d-d_\infty\|^2.
          \label{rate6}
  \eea
Let us define now, for $t\geq 0$,
 \bea y(t)&= & \frac12\|v(t)\|^2+\frac12\|\nabla d(t)-\nabla
  d_\infty\|^2+ \frac{\alpha}{2} \|d(t)-d_\infty\|^2 +\int_\Omega F(d(t))dx-\int_\Omega F(d_\infty) dx\non\\
  &&\ \ -
  \int_\Omega f(d_\infty)(d(t)-d_\infty)dx. \label{y1}
  \eea
In \eqref{ra5} and \eqref{y1}, we choose
  $$ \alpha\geq  1+2C_2 >0,\quad \varepsilon_1=\frac{\nu}{4\alpha}. $$
As a result,
  \be y(t)+C_2\|d-d_\infty\|^2\geq \frac12(\|v\|^2+\|d-d_\infty\|_{H^1}^2).
 \label{y1a}
  \ee
Furthermore, we infer from \eqref{y1a} that for certain constants
$C_3,C_4>0$,
  \be \frac{d}{dt} y(t)+C_3y(t)\leq C_4\|d-d_\infty\|^2\leq
  C(1+t)^{-\frac{2\theta}{1-2\theta}}.\label{ra7}\ee
  As in  \cite{W07,WGZ1}, we have
  \be y(t)\leq C(1+t)^{-\frac{2\theta}{1-2\theta}}, \quad \forall\ t\geq 0,\ee
  which together with \eqref{y1a} implies that
  \be \|v(t)\|+\|d(t)-d_\infty\|_{H^1}\leq C(1+t)^{-\frac{\theta}{1-2\theta}}, \quad \forall\ t\geq 0.\label{rate2}
  \ee
\noindent \textbf{Step 3.} In the last step, we proceed prove the
convergence rate in higher order norm. In Section 3.1, it has been
proven that, once we could obtain the uniform bound of $d$ in $H^2$,
we are able to obtain strong convergence of $d$ in $H^2$ instead of
weak convergence. By reinvestigating the higher order energy
estimate for the subtracted system \eqref{11}--\eqref{33} (cf. also
Lemma \ref{he2d}), we can obtain a further result, which provides
the same rate estimate of $(v,d)$ in $H^1\times H^2$ as
\eqref{rate2}.

In what follows, we just perform the estimates for classical
solutions. Taking the time derivative of $A(t)$, we obtain by a
direct calculation
 \bea
 && \frac12\frac{d}{dt}A(t)+ (\|S v\|^2+ \|\nabla(\Delta
 d-f(d)\|^2)\non\\
 &=&  (S v, v\cdot\nabla v) -(f'(d)(\Delta d-f(d)),
 \Delta d-f(d))+ 2\int_\Omega (\Delta d-f(d))_{x_j}v_{x_k}^jd_{x_k}
 dx\non\\
 && \ \ -\int_\Omega \nabla \pi\cdot \nabla
 d(\Delta d-f(d)) dx\non\\
 &=&(S v, v\cdot\nabla v) -(f'(d)(\Delta d-f(d)),
 \Delta d-f(d))\non\\&&+ 2\int_\Omega (\Delta d-f(d))_{x_j}v_{x_k}^j(d-d_\infty)_{x_k}
 dx
 -2 \int_\Omega (\Delta d-f(d))\nabla v\nabla^2 d_\infty dx\non\\
 && \ \ -\int_\Omega \nabla \pi \cdot \nabla
 d(\Delta d-f(d)) dx\non\\
 &:=& I_2+I_3+I_4+I_5+I_6.
 \label{3r1}
 \eea
 In the above, we use the fact that $(S v, v_t)=(-\Delta v, v_t)$,
 which follows from $v_t\in H$. Noticing that we have got uniform bounds for $\|v\|_{H^1}$ and
$\|d\|_{H^2}$ before (see Lemma \ref{vcon}), in what follows we
estimate $I_i\  (i=2,...,6)$ term by term.
 \bea
 |I_2|&\leq& \|S v\|\|v\|_{L^4}\|\nabla  v\|_{L^4}\leq  C\|S v\|(\|\nabla v\|^\frac12\|v\|^\frac12)(\|\Delta
 v\|^\frac12\|\nabla v\|^\frac12)\non\\
 &\leq& C\|S v\|\|\Delta v\|^\frac12\|v\|^\frac12 \leq \varepsilon_2 \|S v\|^2+C\|v\|^2.\label{3r2}
 \eea
 Since
 \bea
 && \|\nabla (\Delta d-\Delta d_\infty)\|\non\\
  &\leq &\|\nabla (\Delta
 d-f(d))\| +\|\nabla(f(d)-f(d_\infty))\|\non\\
 &\leq& \|\nabla (\Delta
 d-f(d))\|+ \|f'(d)(\nabla d-\nabla d_\infty)\|+\|(f'(d)-f'(d_\infty))\nabla
 d_\infty\|\non\\
 &\leq& \|\nabla (\Delta
 d-f(d))\|+ \|f'(d)\|_{L^\infty}\|(\nabla d-\nabla d_\infty)\|+ \|f''(\xi)\|_{L^\infty}\|d-d_\infty\|_{L^4}\|\nabla
d_\infty\|_{L^4}\non\\
&\leq & \|\nabla (\Delta
 d-f(d))\|+ C\|d- d_\infty\|_{H^1},\label{3r3a}
 \eea
 we have
 \bea
 |I_3|&\leq & \|f'(d)\|_{L^\infty}\|\Delta d-f(d)\|^2\leq C(\|\Delta
 d-\Delta d_\infty\|^2+\|f(d)-f(d_\infty)\|^2)\non\\
 &\leq& C\|\Delta
 d-\Delta d_\infty\|^2+
 C\|f'(\xi)\|^2_{L^\infty}\|d-d_\infty\|^2\non\\
&\leq& C\|\nabla (\Delta d-\Delta
d_\infty)\|^\frac43\|d-d_\infty\|^\frac23+ C\|d-d_\infty\|^2\non\\
&\leq& \varepsilon_2\|\nabla (\Delta
 d-f(d))\|^2+C
\|d-d_\infty\|^2. \label{3r3}
 \eea
 Next,
 \bea
 |I_4|&\leq & \|\nabla (\Delta d-f(d))\|\|\nabla v\|_{L^4}\|\nabla
 (d-d_\infty)\|_{L^4}\non\\
 &\leq & \varepsilon_2 \|\nabla (\Delta d-f(d))\|^2\non\\&&
 + C(\|\Delta
 v\|\|\nabla v\|+\|\nabla v\|^2)(\|\Delta(d-d_\infty)\|\|\nabla
 (d-d_\infty)\|+\|\nabla (d-d_\infty)\|^2)\non\\
 &\leq & \varepsilon_2 \|\nabla (\Delta d-f(d))\|^2+
 \varepsilon_2\|S
 v\|^2 +C\|\Delta (d-d_\infty)\|^2+C\left(1+\frac{1}{\varepsilon_2}\right)\|\nabla(d-d_\infty)\|^2\non\\
 &\leq& 2\varepsilon_2 \|\nabla (\Delta d-f(d))\|^2+
 \varepsilon_2\|S
 v\|^2 +C\|d-d_\infty\|_{H^1}^2
 .\label{3r4}
 \eea
 \bea
 |I_5|&\leq& \|\Delta d-f(d)\|_{L^4}\|\nabla v\|_{L^4}\|\nabla
 ^2d_\infty\|\non\\
 &\leq& C\|\nabla (\Delta d-f(d))\|(\|\Delta v\|^\frac34\|
 v\|^\frac14+\| v\|)\non\\
 &\leq& \varepsilon_2\|\nabla (\Delta d-f(d))\|^2+
 \frac{C}{\varepsilon_2}(\|\Delta v\|^\frac32\| v\|^\frac12+\|v\|^2)\non\\
 &\leq& \varepsilon_2\|\nabla (\Delta d-f(d))\|^2+
 \varepsilon_2\|S
 v\|^2+C\left(\frac{1}{\varepsilon_2^7}+\frac{1}{\varepsilon_2}
 \right)\| v\|^2.\label{3r5}
 \eea
 \bea
 |I_6|&\leq& \|\nabla \pi\|\|\nabla d\|_{L^4}\|\Delta
 d-f(d)\|_{L^4}\non\\
 &\leq& C \|S v\| \|\nabla (\Delta
 d-f(d))\|^\frac12\|\Delta d-f(d)\|^\frac12\|d\|_{H^2}\non\\
 &\leq & \varepsilon_2\|\nabla (\Delta d-f(d))\|^2+
 \varepsilon_2\|S
 v\|^2+ \frac{C}{\varepsilon_2^3} \|\Delta d-f(d)\|^2\non\\
 &\leq& 2\varepsilon_2\|\nabla (\Delta d-f(d))\|^2+
 \varepsilon_2\|S
 v\|^2+ \frac{C}{\varepsilon_2^7}\|d-d_\infty\|_{H^1}^2.
 \label{3r5a}
 \eea
Taking $\varepsilon_2$ sufficiently small, we deduce from
\eqref{3r1}--\eqref{3r5a} that
 \bea \frac{d}{dt}A(t)+ (\|S v\|^2+ \|\nabla (\Delta
 d-f(d))\|^2)\leq C (\|v\|^2+\|d-d_\infty\|_{H^1}^2).\label{3r6}
  \eea
Using the Poincar\'e inequality for $\Delta d-f(d)$ whose trace on
$\Gamma$ is 0 and Lemma \ref{S}, we can conclude from \eqref{3r6}
and \eqref{rate2}
  that
\bea
 \frac{d}{dt}A(t)+ C A(t)\leq C (
 \|v\|^2+\|d-d_\infty\|_{H^1}^2)\leq C(1+t)^{-\frac{2\theta}{1-2\theta}},
 \quad \forall\ t\geq 0.\label{rateaa}
    \eea
 Again, by the Gronwall inequality, we have
  \be A(t)\leq C(1+t)^{-\frac{2\theta}{1-2\theta}}, \quad \forall\ t\geq 0,\ee
  which yields
  \be \|\nabla v(t)\|+\|\Delta d(t)-f(d(t))\|\leq C(1+t)^{-\frac{\theta}{1-2\theta}}, \quad \forall\ t\geq 0.\label{rate3}
  \ee
Recalling \eqref{kkk}, it follows from \eqref{rate3} that
 \be \|\Delta d(t)-\Delta d_\infty\| \leq C(1+t)^{-\frac{\theta}{1-2\theta}}, \quad \forall\ t\geq 0.\label{rate4}
  \ee
Summing up, from \eqref{rate2}\eqref{rate3}\eqref{rate4} we can
deduce the required estimate \eqref{rate}. The proof of Theorem
\ref{main2d} is complete.

\section{Results for Three Dimensional Case }
\setcounter{equation}{0} The results proved in previous section hold
true for global classical solutions to system \eqref{1}--\eqref{5}
in 3-D case. In what follows, we show the convergence to equilibrium
for two subcases considered in \cite{LL95} (ref. Theorem B and
Theorem C therein) that existence of global classical solution was
proven. In particular, we answer the question of uniqueness of
asymptotic limit of $d$ (cf. \cite[Remark, page 32]{LL95}) and
provide a uniform convergence rate.\\

 \textbf{Case I: Initial Data Near Absolute Minimizer of $E$.} \\
The following result has been proven in \cite[Proposition
5.2]{LL95}.

\begin{proposition}
There is an $\varepsilon_0\in (0,1)$ depending only on
$\nu,\lambda,\gamma,\Omega$ and $f$ with the following property:
Whenever
$$ \nu\|\nabla v\|^2(0)+\lambda\gamma\|\Delta d-f(d)\|^2(0)\leq
\varepsilon_0,$$
 either\\
 (1) Problem \eqref{1}--\eqref{5} has a unique classical
 solution $(v,d)$ in $\Omega\times(0,+\infty)$\\
 or\\
 (2) there is a $T_*\in (0, +\infty)$ such that
 $$E(T_*)<E(0)-\varepsilon_0,$$
 where
 $$ E(t)=\|v\|^2+\lambda\|\nabla d\|^2+2\lambda\int_\Omega F(d)dx.$$
 Moreover, in case (1), one has
 \be \|v(t)\|_{H^1(\Omega)}\to 0,\quad \|\Delta d-f(d)\|\to 0,\quad
 \text{as}\ t\to +\infty.
 \ee
\end{proposition}

\noindent Before proving the convergence result corresponding to
Theorem \ref{main2d}, we turn to the second case. Later we shall
prove our result in a unified way.\\

\textbf{Case II: Arbitrary Initial Data with Large Viscosity. }

It has been proven that for any initial data $v_0\in H^1(\Omega)$,
$d_0\in H^2(\Omega)$, if the viscosity $\nu$ is "large enough" (see
below), problem \eqref{1}--\eqref{5} admits a unique global
classical solution (cf. \cite[Theorem B]{LL95}). As pointed out in
\cite{LL95}, when the dimension is three, the size of viscosity
$\nu$ plays a rather crucial role while the other constants
$\lambda,\gamma$ do not, as long as $\lambda,\gamma$ are positive
constants. Thus we shall assume $\lambda=\gamma=1$ for the sake of
simplicity. The following high order energy estimate can be obtained
(cf. \cite[(4.13)]{LL95}).

\bl \label{he3dd} In the 3-D case, the following inequality holds
for
 classical solution $(v,d)$ to problem \eqref{1}--\eqref{5}
 \be \frac12\frac{d}{dt}\tilde{A}(t)\leq -\left(\nu-K\nu^\frac12\tilde{A}\right)\|\Delta
 v\|^2-\left(1-\frac{K\tilde{A}}{\nu}\right)
 \|\nabla (\Delta d-f(d))\|^2+K\tilde{A}, \quad \forall\
 t\geq 0,\label{he3d}
 \ee
 where $\tilde{A}=A+1=\|\nabla v\|^2+ \|\Delta d-f(d)\|^2+1$ (cf. \eqref{A})
 and $K$ is a positive constant depending on $f,\nu, \Omega, \|v_0\|, \|d_0 \|_{H^1(\Omega)}$.
 \el

When the viscosity $\nu$ is assumed to be properly large, based on
the above lemma, we can not only show that the global solution
$(v,d)$ is uniformly bounded (as in \cite{LL95}) but also the
quantity $A(t)$ decays to zero in time.

It follows from \eqref{ae} that
 \be \int_t^{t+1} \tilde{A}(\tau)d\tau\leq \int_t^{t+1}
 A(\tau)d\tau+ 1\leq M,\quad \forall\ t\geq 0,\label{AM}
 \ee
 where $M>0$ is a constant depending only on $\|v_0\|,
 \|d_0\|_{H^1}$. Then we have

 \bl If \be \nu^\frac12\geq
 K\left(\tilde{A}(0)+2KM+4M\right)+\frac12,\label{nu}\ee
then the unique global solution to problem \eqref{1}--\eqref{5}
satisfies
 the following uniform estimate
  \be \|v(t)\|_{H^1}+
\|d(t)\|_{H^2}\leq C, \quad \forall\  t\geq 0,\label{ubd}
 \ee
 where $C$ is a constant depending on $f, \Omega, \|v_0\|_{H^1(\Omega)}, \|d_0
 \|_{H^2(\Omega)}$.  Furthermore,
 \be \lim_{t\rightarrow +\infty} (\|v(t)\|_{H^1}+ \|-\Delta d(t)+f(d(t))\|)=0. \label{vcon3d}\ee
 \el
\begin{proof}
Proof of existence and uniqueness of the global solution has been
given in \cite{LL95}. Next, we show the uniform bound \eqref{ubd}.
Take $\nu$ large enough that \eqref{nu} is satisfied. Then by
\eqref{he3d}, there must be some $T_0>0$ such that
$$ \nu-K\nu^\frac12\tilde{A}(t)\geq 0,\quad
1-\frac{K\tilde{A}(t)}{\nu}\geq 0,$$  for all $t\in [0,T_0]$.
Moreover, on $[0,T_0]$,
 \be \frac{d}{dt}\tilde{A}(t)\leq 2K\tilde{A}(t).\label{at}
 \ee
 Denote $T_*=\sup T_0$. First we show that $T_*\geq
 1$ by a contradiction argument.\\
 If $T_*<1$, then
 $$\tilde{A}(T_*)\leq \tilde{A}(0)+2K\int_0^1\tilde{A}(t)dt\leq
 \tilde{A}(0)+2KM.$$
  On the other hand, from the definition of $T_*$, we have
 $$ \nu<\max\{K\tilde{A}(T_*), K^2\tilde{A}^2(T_*)\}\leq K(\tilde{A}(0)+2KM)+ K^2(\tilde{A}(0)+2KM)^2,$$
which contradict \eqref{nu}.

Next, if $T_*<+\infty$, \eqref{AM} implies that there is a $t_1\in
[T_*-\frac12,T_*]$ such that
 \be \tilde{A}(t_1)\leq 4M.\ee
 As a result,
 \be \tilde{A}(T_*)\leq 4M+2K\int_{t_1}^{T_*}\tilde{A}(t)dt\leq
 4M+2KM.\label{TT}\ee
Again from the definition of $T_*$, we have
$$ \nu<\max\{K\tilde{A}(T_*), K^2\tilde{A}^2(T_*)\},$$
which together with \eqref{TT} yields  a contradiction with
\eqref{nu}.

Therefore, for all $t\geq 0$, \eqref{at} holds. Namely,
 \be \frac{d}{dt}A(t)\leq 2KA(t)+2K\leq KA^2(t)+ 3K.\ee
 Due to
\eqref{ae}, we can conclude \eqref{ubd} and \eqref{vcon3d} following
the similar argument in the proof of Lemma \ref{vcon}.
\end{proof}

\br Generally speaking, \eqref{nu} only provides a sufficient
condition on the largeness of viscosity $\nu$, which ensures the
existence of global solution to problem \eqref{1}--\eqref{5}. It may
not be an optimal lower bound for all possible $\nu$. \er

Based on above results, now for both cases I and II, one can argue
exactly as in Section 3.1 to conclude
 \be \lim_{t\rightarrow +\infty}
 (\|v(t)\|_{H^1}+\|d(t)-d_\infty\|_{H^2})=0.\label{cgce3d}
 \ee

 Then we are able to proceed to show the estimate on
convergence rate for both two cases. To this aim, we check the
argument for 2-D case step by step. By applying corresponding
Sobolev embedding Theorems in 3-D, we can see that all calculations
in Section 3.2 are valid for our current case (with minor
modifications). Hence the details are omitted.

We complete the proof for Theorem \ref{main3da}  and Theorem
\ref{main3d}.

\section{Further Remarks}\setcounter{equation}{0} We remark that our approach used in this paper are valid for some
other model systems for nematic liquid crystal flows in the
literature and similar convergence result can be proved.

\noindent (1) \textbf{A model with changing density}

Recently the following problem was considered in \cite{LZ08}.
 \bea \rho_t+\nabla\cdot(\rho v)&=&0,\quad \rho\geq 0,\label{LZ1}\\
 (\rho v)_t+\nabla\cdot(\rho v\odot v)-\nu \Delta v+\nabla P&=&-\lambda
 \nabla\cdot(\nabla d\odot\nabla d),\\
 \nabla \cdot v &=& 0,\\
 d_t+v\cdot\nabla d&=&\gamma(\Delta d-f(d)),\label{LZ2}
 \eea
in $\Omega \times(0,\infty)$, where $\Omega \subset \mathbb{R}^n
(n=2,3)$ is a bounded domain with smooth boundary $\Gamma$.
$\rho(x,t)$ is a scalar function denoting the density of the fluid.
The above density-dependent liquid crystal model is subject to the
following initial condition
 \be
 \rho|_{t=0}=\rho_0(x)\geq 0,\ \ (\rho v)|_{t=0}=q_0(x),\quad
 d|_{t=0}=d_0(x),\qquad \text{for}\ x\in \Omega,\label{LZ3}
 \ee
 and the boundary conditions:
 \be
 v(x,t)=0,\quad d(x,t)=d_0(x),\qquad \text{for}\ (x, t)\in \Gamma\times
 \mathbb{R}^+.
 \label{LZ4}
 \ee

 Problem \eqref{LZ1}--\eqref{LZ4} can be viewed as a generalization of our problem
 \eqref{1}--\eqref{5}. It enjoys some important properties as for \eqref{1}--\eqref{5}.
In particular, we have the following \emph{basic energy law} (see
\cite{LZ08})
 \be \frac{d}{dt}\int_\Omega \left(\frac{1}{2}\rho |v|^2+\frac{\lambda}{2}|\nabla
 d|^2+\frac{\lambda}{2}\int_\Omega F(d)dx\right)=-\int_\Omega
 \left(\nu |\nabla
 v|^2+\lambda\gamma|\Delta d-f(d)|^2\right)dx.\label{LZe}
 \ee

 In \cite{LZ08}, the authors proved the
existence of the weak solution to incompressible liquid crystal
system \eqref{LZ1}--\eqref{LZ4} under certain compatibility
condition on the initial data. There they considered the general
case for the density, namely they only required the initial density
to be nonnegative. As a result, one can only expect the density to
be nonnegative for all time and vacuum state may occur. In this
case, it is very difficult to prove corresponding results to
\cite[Theorem B, Theorem C]{LL95} where the global existence and
uniqueness as well as asymptotic behavior of classical solutions
were obtained.

However, if we assume in addition that the initial density is a
bounded positive function, i.e., there are two positive constants
$\underline{\rho}$ and $\bar{\rho}$ such that
 \be 0<\underline{\rho}\leq \rho_0(x)\leq \bar{\rho},\quad \forall\
 x\in \Omega.\ee
 Then by virtue of the comparison principle (cf. \cite{LZ08}), we
 have
 \be 0<\underline{\rho}\leq \rho(x,t)\leq \bar{\rho},\quad \forall\
 t\geq 0.\label{rho}
 \ee

 In this special case, we can check that under suitable assumptions
 on the initial data (for instance, $v_0\in H_0^1$, $d_0\in H^2$),
 parallel results to \cite[Theorem B,
Theorem C]{LL95} can be achieved. Besides the basic energy law
\eqref{LZe}, due to uniform upper and lower bound \eqref{rho}, one
can proceed to get proper high order energy law similar to Lemma
\ref{he2d} (2-D case) as well as Lemma \ref{he3dd} (3-D case). For
instance, denote
 \be
 \hat{A}(t)=\rho\|\nabla v\|^2+ \lambda\|\Delta d-f(d)\|^2,\label{A1}
 \ee
 we can show that
 \bl In the 2-D case, the following inequality holds for the
 classical solution $(v,d)$ to problem \eqref{LZ1}--\eqref{LZ4}
  \be \frac{d}{dt}\hat{A}(t) +K_1(\|S v\|^2+\|\nabla(\Delta d-f(d))\|^2) \leq K_2(\hat{A}^2(t)+1), \qquad \forall
 t\geq 0,\label{hee}
 \ee
 where $K_1,K_2$ are constants depending on
  $f, \Omega, \|v_0\|, \|d_0 \|_{H^1(\Omega)}, \nu,\lambda, \gamma, \bar{\rho},\underline{\rho}$.
 \el

 Corresponding results in 3-D case (cf. Lemma \ref{he3dd}) can also be obtained. The proofs for these results
 follow from the same sprit of those in \cite{LL95} with some proper
 modifications. Hence, the details are omitted
 here.
 \br
Different from system \eqref{1}--\eqref{5}, where the density is
assumed to be a constant, in order to get the high order energy law,
we deal with the time derivative of a modified quantity $\hat{A}(t)$
with weight $\rho$ instead of $A(t)$. This is due to the
mathematical structure of \eqref{LZ1}--\eqref{LZ4}, in which the
density variable is involved. Because of the uniform upper and lower
bounds of the density \eqref{rho}, one can check that $\hat{A}(t)$
plays a similar role as $A(t)$ for system \eqref{1}--\eqref{5}.
 \er

 Based on the facts obtained above, we are able to prove the
 corresponding convergence results (cf. Theorem 1.1--Theorem 1.3) for system
 \eqref{LZ1}--\eqref{LZ4},
 following the argument in the previous sections. We leave the details to the interested readers.\\

\noindent (2) \textbf{A model with free-slip boundary condition}
 \bea
 v_t+v\cdot\nabla v-\nu  \texttt{div}D(v)+\nabla P&=&-\lambda
 \nabla\cdot(\nabla d\odot\nabla d),\label{1la}\\
 \nabla \cdot v &=& 0,\label{2la}\\
 d_t+v\cdot\nabla d&=&\gamma(\Delta d-f(d)),\label{3la}
 \eea
in $\Omega \times(0,\infty)$, where $\Omega \subset \mathbb{R}^n
(n=2,3)$ is a bounded polygonal domain (with piecewise smooth
boundary). $D(v) = \frac12 (\nabla v + (\nabla v)^T)$ is the
stretching tensor. We consider the system \eqref{1la}--\eqref{3la}
subject to the initial conditions
 \be
 v|_{t=0}=v_0(x) \ \ \text{with}\ \nabla\cdot v_0=0,\quad
 d|_{t=0}=d_0(x),\qquad \text{for}\ x\in \Omega,\label{4la}
 \ee
 and the free-slip boundary conditions:
 \be
 v\cdot \textbf{n}=0,\quad (\nabla\times v)\times \textbf{n}=0,\quad
 \partial_\textbf{n} d=0,\qquad \text{for}\ (x, t)\in \Gamma\times \mathbb{R}^+,
 \label{5la}
 \ee
where $\textbf{n}$ is the unit outer normal vector to the boundary
$\Gamma$.

As has been pointed out in the recent paper \cite{LS01}, the
free-slip boundary condition \eqref{5la} indicates that in the
liquid crystal flows, there is no contribution from the director
field $d$ to the surface forces. Boundary condition \eqref{5la}
seems to be more appropriate for some types of flow in the bulk of a
liquid crystal configuration. On the other hand, it allows people to
construct more efficient numerical schemes for the numerical
simulations for liquid crystal flows (cf. \cite{LS01}). Comparing
with system \eqref{1}--\eqref{5}, the influences of the corner
singularities is less severe with free-slip and Neumann boundary
conditions than the Dirichlet boundary conditions.

Basic theoretical analysis on problem \eqref{1la}--\eqref{5la} has
been done in \cite{LS01}, where the authors proved global existences
of weak solutions as well as regularities and global
existence/uniqueness of classical solutions. In particular, although
the boundary condition \eqref{5la} plays a significantly different
role in the calculation, proper high order energy law similar to
Lemma \ref{he2d} could still be obtained (cf. \cite[Lemma
4.1]{LS01}). The same convergence results for system
\eqref{1}--\eqref{5} obtained in the present paper can be shown true
for problem \eqref{1la}--\eqref{5la}, by adapting the argument here.
We thus omit the details.

\medskip

{\bf Acknowledgement.}  A part of this paper was written during H.
Wu's visit to the Department of Mathematics in Penn State
University. The hospitality of Prof. C. Liu and the department is
gratefully acknowledged. The research of H. Wu was partially
supported by China Postdoctoral Science Foundation.

\end{document}